\documentclass[11pt]{article}
\usepackage{amsmath,amssymb}
\usepackage[dvips]{graphics}
\title{\Large\bf 
Fundamental solution and the weight functions of
the transient problem on a semi-infinite crack propagating in a half-plane}
\author{ Y.A. Antipov, A.V. Smirnov\\ 
Department of Mathematics, Louisiana State University\\
Baton Rouge LA 70803\\
 }

\date{}

\newcommand{\sgn}{\mathop{\rm sgn}\nolimits}

\newcommand{\I}{\mathop{\rm Im}\nolimits}
\newcommand{\R}{\mathop{\rm Re}\nolimits}
\newcommand{\const}{\mbox{const}}

\newcommand{\beqa}{\begin{eqnarray}}
\newcommand{\eeqa}[1]{\label{#1}\end{eqnarray}}
\newcommand{\bequ}{\begin{equation}}
\newcommand{\eequ}[1]{\label{#1}\end{equation}}

\newcommand{\Md}{\partial}

\newcommand{\Ga}{\alpha}
\newcommand{\Gb}{\beta}
\newcommand{\Gd}{\delta}

\newcommand{\Gve}{\varepsilon}

\newcommand{\Gvf}{\varphi}
\newcommand{\Gg}{\gamma}
\newcommand{\Gc}{\chi}

\newcommand{\Gk}{\kappa}

\newcommand{\Gl}{\lambda}

\newcommand{\Gt}{\theta}

\newcommand{\Gr}{\rho}

\newcommand{\Gs}{\sigma}

\newcommand{\Go}{\omega}
\newcommand{\Gx}{\xi}
\newcommand{\Gy}{\psi}
\newcommand{\Gz}{\zeta}
\newcommand{\GD}{\Delta}

\newcommand{\GG}{\Gamma}

\newcommand{\GP}{\Pi}

\newcommand{\GO}{\Omega}

\newcommand{\GY}{\Psi}

\newcommand{\CK}{{\cal K}}
\newcommand{\CL}{{\cal L}}

\newcommand{\CQ}{{\cal Q}}

\newcommand{\CX}{{\cal X}}

\newcommand{\beq}{\begin{equation}}
\newcommand{\eeq}{\end{equation}}
\newcommand{\barr}{\begin{eqnarray}}
\newcommand{\earr}{\end{eqnarray}}
\newcommand{\beqn}{\begin{equation*}}
\newcommand{\eeqn}{\end{equation*}}
\newcommand{\barrn}{\begin{eqnarray*}}
\newcommand{\earrn}{\end{eqnarray*}}

\newcommand{\fr}{\frac}

\begin{document}
\maketitle

\noindent

\begin{abstract}

The two-dimensional  transient  problem that is studied concerns a semi-infinite crack in an 
isotropic solid comprising 
an infinite strip and a half-plane joined together and having  the same elastic constants. The crack propagates along the interface
at constant speed subject to time-independent loading.  
By means of the Laplace and Fourier transforms the problem
is formulated as a vector Riemann-Hilbert problem. When the distance from the crack to the boundary
grows to infinity the problem admits a closed-form solution. In the general case, a method of partial matrix factorization
is proposed. In addition to factorizing some scalar functions it requires solving a certain system
of integral equations whose numerical solution is computed by the collocation method. The stress intensity 
factors and the associated weight functions are derived. Numerical results for
 the weight functions are reported and the boundary effects are discussed.
The weight functions are employed to describe propagation of a semi-infinite crack beneath the half-plane boundary at piecewise constant speed.

\end{abstract}

\setcounter{equation}{0}

\section{Introduction}

The main goal of this work is to propose a new method for analyzing plane dynamic transient problems when the two modes are coupled, and the standard Wiener-Hopf method does not work. In addition
to factorization of two scalar functions it employs derivation and solution of a certain  system of two
integral equations. 
The method is illustrated by the study of a crack propagating at sub-Rayleigh speed parallel to the boundary of a solid when loading 
is time independent. 
The model problem admits formulation as a vector Riemann-Hilbert problem (RHP). 
In the case when the crack is far away from the boundary, the problem can be modeled as propagation
of a semi-infinite crack along the interface between two weakly bonded, identical and isotropic half-planes. The problem on  crack growth in a plane
at constant sub-Rayleigh speed  due to general time-independent  loading (including the case of concentrated forces  applied to the crack faces) 
was solved exactly  by Freund (1990) by means of the Wiener-Hopf method. 
Huang and Gao (2001, 2002)
analyzed the intersonic regime including the case of concentrated forces (the fundamental solution)
and the model problem on a suddenly stopping crack.  
When the crack is close to the boundary of the body, the boundary effects cannot be ignored, and the problem on a crack propagating
parallel to the half-plane boundary can be considered as an adequate model.
In the static case, the matrix coefficient
of the RHP  admits a closed-form  factorization (Zlatin and Khrapkov, 1986). The steady-state case, 
when  loads applied to the propagating crack   move with the crack at the same constant speed,  was recently
analyzed by Antipov and Smirnov (2013). By means of the Fourier transform, the problem was  mapped into a vector RHP whose matrix coefficient
 did not allow for an explicit factorization. 
The RHP was rewritten as a system of singular integral equations, and an approximate method
of orthogonal polynomials for its solution was proposed. To the authors' knowledge,
an analytic solution  to the transient problem on a semi-infinite crack propagating along the boundary of a half-plane is not available in the literature. 

In Section 2, we describe the transient model and apply the Fourier and Laplace transforms in a standard manner
(Freund, 1990; Slepyan, 2002; Antipov and Willis, 2003) in order to reduce the governing boundary-value problem to an order-2 vector RHP.  Although the matrix coefficient has the same structure as in the steady-state
case  (Antipov and Smirnov, 2013), the transient problem is much harder  for the parameters $\Ga$ and $\Gb$ involved
being functions of the Laplace and Fourier parameters not constants as in the steady-state problem.
All our efforts availed us no results not only in factorizing the
matrix coefficient of the RHP, but even in computing the partial indices of factorization (Vekua, 1967).
The partial indices play an essential part in solvability theory of a vector RHP and in theory of approximate Wiener-Hopf matrix factorization. According to the stability criterion 
for partial indices (Bojarski, 1958; Gohberg and Kre\u in, 1958;
Vekua, 1967) applied to a $2\times 2$ matrix the partial indices, integers $\Gk_1$ and $\Gk_2$, are stable if and only if $|\Gk_1-\Gk_2|\le 1$. If they do not satisfy this criterion, then  approximate canonical
Wiener-Hopf factors may not converge to the exact ones.  At the same time, without
knowledge of exact factorization, in general, there is no way to determine the partial indices.
An example (not inspired by an applied physical problem) of unstable partial indices is given by Litvinchuk and Spitkovski\u i (1987).
It turns out that the partial indices associated with  contact, fracture, and diffraction models available in the literature (Moiseyev and Popov, 1990; Antipov and Moiseyev, 1991; Antipov and Silvestrov, 2002)
are stable.   Although this circumstance makes the determination of the partial indices of factorization not
an absolute necessity but rather a desideratum, in this paper we surmount the deficiency of not knowing the partial indices by bypassing the 
problem of approximate matrix factorization. Instead, we propose the method
of partial factorization that comprises factorization of some scalar
functions and numerical solution of a certain system of integral equations

To introduce the reader to the method
presented later in the paper, in Section 3, we analyze  the transient problem for a semi-infinite  crack in the whole plane.
In this case, the vector RHP is decoupled and solved by quadratures. We also derive exact formulas for the stress intensity factors (SIFs)
and the weight functions introduced by
Bueckner (1970)  
for a semi-infinite static
crack in a homogeneous elastic medium. For the elastic case exact and 
approximate expressions for the weight functions are available in the literature
for a variety of models. Exact weight functions for a static semi-infinite
interfacial crack in a three-dimensional
unbounded body  were constructed by Antipov (1999) and later on employing a different approach 
by Bercial-Velez et al (2005). The matrix factorization technique was applied 
by Antipov (2012) for recovering exact representations of the weight functions in
  micropolar theory.  On using a scalar factorization method Willis and Movchan (1995) 
  derived the dynamic Mode I weight function for a semi-infinite crack in a three-dimensional homogeneous body  under the conditions of steady-state normal loading. The case of shear loading
  analyzed by them later   (Movchan and Willis, 1995) revealed coupling between Modes II and III.
 On applying the matrix factorization method they derived exact formulas for the  
 steady-state Mode II and III weight functions.
Mode I and II weight functions for a viscoelastic medium
with different bulk and shear relaxation were determined in the
steady-state case by Antipov and Willis (2007). 
Recently, Antipov and Smirnov (2013) found approximate formulas
for weight functions associated with the steady-state semi-infinite crack propagating 
bellow the boundary of a half-plane. 
The transient weight functions are available  for a semi-infinite crack
propagating at constant speed only for an unbounded body
 in the elastic case for sub-Rayleigh speed (Freund, 1990),
  the transonic regime 
(Huang and Gao, 2001), and in the viscoelastic case when
the constant speed of crack propagation may take any value up to the speed of
dilatational waves (Antipov and Willis, 2003). 
In Section 3 of this paper, we derive  
the transient weight functions for a plane by a factorization method
that will be later generalized for a half-plane.

In Section 4,  we propose an approximate method for the vector RHP associated with 
the transient problem for a half-plane.  First, we
split the matrix coefficient into a diagonal matrix that is discontinuous at infinity and a matrix  that is continuous.
On  factorizing the discontinuous part and recasting the vector RHP we derive a new vector RHP that is equivalent
to a system of two integral equations on the interval $(-\infty,0)$. The diagonal elements of the matrix kernel are constants, while the off-diagonal elements
are continuous functions which have an order-2 zero at infinity. We show that in order to determine the Laplace  transforms of the SIFs and the weight functions, 
it is sufficient  to know the solution  to the system of integral equations at the point 0 only. We describe the numerical procedure and the inversion method 
of the Laplace transform we applied and discuss the numerical results for the weight functions.

With the fundamental solution and weight functions being derived we proceed, in 
Section 5, to describe nonuniform growth of a semi-infinite crack parallel to the
boundary of a half-plane. For the whole plane such an algorithm based on the fundamental solution and the solution of the model problem on a suddenly stopped crack is known
(Freund, 1990).  We aim to generalize this procedure for a half-plane. The main feature here
is to  take into account 
the fact that the Mode I and II weight functions after
the longitudinal wave reflects from the boundary and strikes the crack 
do not act alone anymore and the off-diagonal weight functions play a substantial part. We show that in order to determine the stresses radiated out by a suddenly stopped crack, one needs to solve a system of two Volterra convolution equations,
 a generalization of the single Abel equation appeared in the Freund method for the whole plane. On solving this system exactly we determine
the stresses the crack needs to negate on the prospective fracture plane to proceed further. The procedure to be exposed allows for the possibility of finding the SIFs
at the tip of a crack propagating at piecewise constant speed bellow the boundary.

\setcounter{equation}{0}

\section{Transient problem for a half-plane as a vector RHP}

\subsection{Formulation}

\begin{figure}[t]
\centerline{
\scalebox{0.6}{\includegraphics{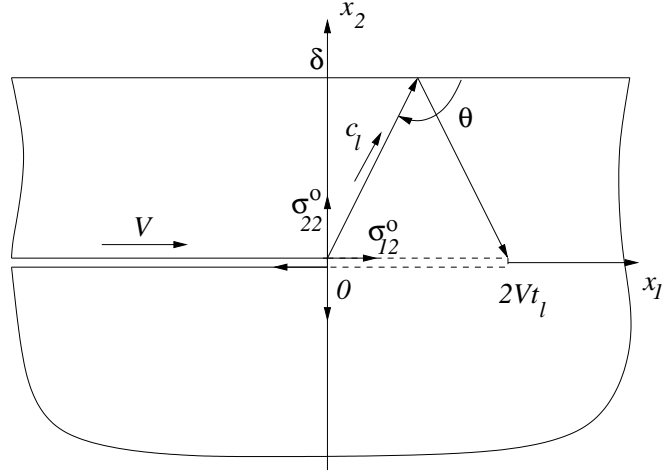}}}
\caption{A semi-infinite crack propagating parallel to the boundary: the transient problem.}
\label{fig0}
\end{figure}

The elastic medium $\GP=\{|x_1|<\infty,\; -\infty<x_2<\Gd\}$ through which the crack propagates  comprises an infinite strip
$\{|x_1|<\infty,\; 0<x_2<\Gd\}$ and a half-plane $\{|x_1|<\infty,\; -\infty<x_2<0\}$ bonded together. The bonding is not perfect, and it is assumed that along the interface there is a semi-infinite crack. The faces of the crack are subjected to plane strain loading that forces the crack to propagate at constant sub-Rayleigh
speed $V$. The presence of the weak interface encourages the crack to propagate parallel
to the boundary $\{|x_1|<\infty, x_2=\Gd\}$ rather than deviate towards it (Fig. \ref{fig0}).
The boundary of the body $\GP$ is free of traction,
and the Lam\'e constants $\lambda$ and $\mu$ and 
the density $\rho$ of the strip and the half-plane
are assumed to be the same. 
  Boundary conditions of the problem are specified to the form
\beq
\begin{aligned}
	&\Gs_{j2}=-\Gs_{j2}^\circ(x_1,0)H(t), &\; -\infty<x_1<Vt,&\quad x_2=0^\pm,\\
	&\Gs_{j2}=0,&\; -\infty<x_1<\infty,&\quad x_2=\delta,
\end{aligned}\quad j=1,2,
\label{1.0}
\eeq
where $\Gs_{12}^\circ$, $\Gs_{22}^\circ$ are prescribed functions and $H(t)$ is the unit step function. 

It is advantageous to change variables from the material coordinates $x_1$, $x_2$ to the crack tip coordinates $x=x_1-Vt$, $y=x_2$. In these coordinates, displacement potentials $\Gvf$ and $\psi$ of the medium satisfy the wave equations
\beq
\begin{aligned}
	c_l^2\hat\Ga^2\fr{\Md^2\Gvf}{\Md x^2}+c_l^2\fr{\Md^2\Gvf}{\Md y^2}
	+2V\fr{\Md^2\Gvf}{\Md x\Md t}-\fr{\Md^2\Gvf}{\Md t^2}=0,\\
	c_s^2\hat\Gb^2\fr{\Md^2\psi}{\Md x^2}+c_s^2\fr{\Md^2\psi}{\Md y^2}
	+2V\fr{\Md^2\psi}{\Md x\Md t}-\fr{\Md^2\psi}{\Md t^2}=0, 
	\end{aligned}
	\quad (x,y)\in\hat\GP\setminus\{-\infty<x<0, y=0\}, \quad t>0,
	\label{1}
\eeq
where $\hat\GP=\{|x|<\infty, -\infty<y<\Gd\}$, subject to the zero initial conditions
\beq
\Gvf=\psi=0,\quad \fr{\Md\Gvf}{\Md t}=\fr{\Md\psi}{\Md t}=0, \quad (x,y)\in\hat \GP, \quad t< 0.
\label{2}
\eeq
Here, $c_l$ and $c_s$ are the longitudinal and shear wave speeds
\beq
c_l=\sqrt{\fr{\Gl+2\mu}{\Gr}}, \quad 
c_s=\sqrt{\fr{\mu}{\Gr}},
\label{3}
\eeq
and
\beq
\hat\Ga=\sqrt{1-v_l^2}, \quad \hat\Gb=\sqrt{1-v_s^2}, \quad v_l=V/c_l, \quad v_s=V/c_s.
\label{4}
\eeq
In the new coordinates, the  displacements  $u$ and $v$ and the stresses $\Gs_{xy}$ and $\Gs_{yy}$
are expressed through the
dynamic potentials 
$\Gvf$ and $\psi$ as
$$
u=\fr{\Md\Gvf}{\Md x}+\fr{\Md\Gy}{\Md y},\quad v=\fr{\Md\Gvf}{\Md y}-\fr{\Md\Gy}{\Md x},
$$
$$
\Gs_{xy}=\mu\left(\fr{2\Md^2\Gvf}{\Md x \Md y}-\fr{\Md^2\Gy}{\Md x^2}+\fr{\Md^2\Gy}{\Md y^2}\right),
$$
\beq
\Gs_{yy}=\Gl\fr{\Md^2\Gvf}{\Md x^2}+(\Gl+2\mu)\fr{\Md^2\Gvf}{\Md y^2}-2\mu\fr{\Md^2\Gy}{\Md x\Md y}.
\label{5}
\eeq
To complete the formulation of the problem, we explicitly write the boundary conditions in the moving coordinates
$$
\Gs_{xy}=0,\quad \Gs_{yy}=0, \quad |x|<\infty, \quad y=\Gd, \quad t\ge 0,
$$
\beq
\Gs_{xy}=-\Gs_1^\circ(x+Vt,0) \quad \Gs_{yy}=-\Gs_2^\circ(x+Vt,0), 
\quad  -\infty<x<0,  \quad y=0^\pm, \quad t\ge 0.
\label{6}
\eeq
Notice that in the new coordinates the loading is time-dependent.

\subsection{Vector RHP}

We next aim to transform the boundary value problem
to a vector RHP. On 
applying first the Laplace transform with respect to time
\beq
	\left(
	\begin{array}{c}
	\hat\Gvf\\
	 \hat\psi\\
	 \end{array}\right)(x,y,s)=
	\int_0^\infty 
	\left(
	\begin{array}{c}
	\Gvf\\
	 \psi\\
	 \end{array}\right)
	(x,y,t)e^{-s t}dt,\quad \R s=\Gs>0,
	\label{7}
\eeq
and then the Fourier transform with respect to $x$
\beq
\left(
	\begin{array}{c}
	\tilde\Gvf\\
	 \tilde\psi\\
	 \end{array}\right)	
	(p,y,s)=
	\int_{-\infty}^\infty 
	\left(
	\begin{array}{c}
	\hat\Gvf\\
	 \hat\psi\\
	 \end{array}\right)
	(x,y,s)e^{ipx}dx,\quad p\in{\mathbb R},
	\label{8}
\eeq
we can write the governing equations in the form
\beq
	\frac{\Md^2\tilde\Gvf}{\Md y^2}-\alpha^2\tilde\Gvf=0,\quad
	\frac{\Md^2\tilde\psi}{\Md y^2}-\beta^2\tilde\psi=0,\quad
	y\in\{-\infty,\delta\}\setminus\{0\},
	\label{9}
\eeq
where 
\beq
\alpha^2=\hat\Ga^2p^2+2ips v_l/c_l+s^2/c_l^2,\quad
\beta^2=\hat\Gb^2p^2+2ips v_s/c_s+s^2/c_s^2.
\label{10}
\eeq
To fix single branches of the two-valued functions (\ref{10}), we cut the $p$-plane
along lines that pass through the infinite point and join the branch
points $a_\pm$ of the former function and $\Gb_\pm$ of the second one, 
\beq
a_\pm=\fr{is}{V\pm c_l}\in{\mathbb C^\pm}, \quad b_\pm=\fr{is}{V\pm c_s}\in{\mathbb C^\pm}.
	\label{11}
\eeq
We denote the single branches as
\beq
\Ga=\hat\Ga(p-a_-)^{1/2}(p-a_+)^{1/2}, \quad
\Gb=\hat\Gb(p-b_-)^{1/2}(p-b_+)^{1/2},
\label{12}
\eeq
assuming that as $p\in L=(-\infty,+\infty)$,
$$
-\pi<\arg(p-a_+)<0,\quad 0<\arg(p-a_-)<\pi,
$$
\beq
-\pi<\arg(p-b_+)<0,\quad 0<\arg(p-b_-)<\pi.
	\label{13}
\eeq
In these circumstances, (\ref{12}) implies $\R\Ga>0$ and $\R\Gb>0$ as $p\in L$.
Then the general solution to the differential equations (\ref{9}), bounded as $y\to-\infty$, reads
\beq
	\tilde\Gvf(p,y,s)=C_0(p,s)e^{\alpha y},\quad \tilde\psi(p,y,s)=D_0(p,s)e^{\beta y},\quad -\infty<y<0,
	\label{15}
\eeq
and
\beq
	\begin{aligned}
		\tilde\Gvf(p,y,s)=C_1(p,s)\cosh(\alpha y)+C_2(p,s)\sinh(\alpha y),\\
		\tilde\psi(p,y,s)=D_1(p,s)\cosh(\beta y)+D_2(p,s)\sinh(\beta y),
	\end{aligned}\quad 0<y<\delta.
	\label{16}
\eeq
It is helpful to introduce new functions representing the jumps of the tangential derivatives of the displacement components $u$, $v$ on the crack faces 
\beq
	\begin{aligned}
		&\frac{\Md u}{\Md x}(x,0^+,t)-\frac{\Md u}{\Md x}(x,0^-,t)=\chi_1(x,t),\\
		&\frac{\Md v}{\Md x}(x,0^+,t)-\frac{\Md v}{\Md x}(x,0^-,t)=\chi_2(x,t),
	\end{aligned}\quad -\infty<x<0,\quad t>0.
	\label{17}
\eeq
Then we define the Laplace transforms with respect to time
$$
\hat\Gc_j(x,s)=\int_0^\infty\Gc_j(x,t)e^{-st}dt,\quad
q_j(x,s)=\int_0^\infty\Gs^\circ_j(x+Vt,0)e^{-st}dt,\quad x<0,\quad j=1,2,
$$
\beq
\hat\Gs_1(x,s)=\int_0^\infty\Gs_{xy}(x,0,t)e^{-st}dt,\quad
\hat\Gs_2(x,s)=\int_0^\infty\Gs_{yy}(x,0,t)e^{-st}dt,\quad x>0,
\label{18}
\eeq
and the one-sided Fourier transforms
$$
\tilde\Gc^-_j(p,s)=\int_{-\infty}^0\hat\Gc_j(x,s)e^{ipx}dx,\quad
\tilde q_j^-(p,s)=\int_{-\infty}^0 q_j(x,s)e^{ipx}dx,
$$
\beq
\tilde\Gs^+_j(p,s)=\int_{0}^\infty\hat\Gs_j(x,s)e^{ipx}dx,\quad j=1,2.
\label{19}
\eeq
In order to derive the governing vector RHP, we apply the Laplace and Fourier transforms  
to the six boundary conditions (\ref{6}) and equations (\ref{17}), use the notations
(\ref{18}) and (\ref{19}), and eliminate the functions $C_j(p,s)$ and $D_j(p,s)$ ($j=0,1,2$). The two equations left
comprise the vector RHP
\beq
	\left(\begin{array}{c}\tilde\sigma_1^+(p,s)\\\tilde\sigma_2^+(p,s)\end{array}\right)=\mu i G(p,s)\left(\begin{array}{c}\tilde\chi^-_1(p,s)\\\tilde\chi_2^-(p,s)\end{array}\right)
+\left(\begin{array}{c}\tilde q_1^-(p,s)\\\tilde q_2^-(p,s)\end{array}\right),\quad p\in L.
	\label{20}
\eeq
The functions $\tilde\sigma_j^+(p,s)$ 
are analytic in the upper $p$-half-plane and vanish as 
$p\to\infty$, while the functions $\tilde\chi_j^-(p,s)$ and $\tilde q_j^-(p,s)$ are analytic in the lower $p$-half-plane  and also vanish at infinity. The matrix coefficient of the problem is defined by
\beq
	G(p,s)=\left(\begin{array}{cc}
		g_{11}(p,s) & i g_{12}(p,s)\\
		-i g_{12}(p,s) & g_{22}(p,s)
	\end{array}\right),
	\label{21}
\eeq
$$\begin{aligned}
	&g_{11}(p,s)=\frac{e^{-(\alpha+\beta)\delta}}{2\beta(p^2-\beta^2)p}\left[R_1\sinh\{(\alpha+\beta)\delta\} - R_2\sinh\{(\alpha-\beta)\delta\} + \frac{2\Delta}{R_1}\right],\\
	&g_{12}(p,s)=\frac{4R_2(p^2+\beta^2)}{R_1(p^2-\beta^2)}e^{-(\alpha+\beta)\delta}\sinh^2\frac{(\alpha-\beta)\delta}2,\\
	&g_{22}(p,s)=\frac{e^{-(\alpha+\beta)\delta}}{2\alpha(p^2-\beta^2)p}\left[R_1\sinh\{(\alpha+\beta)\delta\} + R_2\sinh\{(\alpha-\beta)\delta\} + \frac{2\Delta}{R_1}\right],
\end{aligned}$$
\beq
	\begin{aligned}
		&\Delta=R_1^2\sinh^2\frac{(\alpha+\beta)\delta}2-R_2^2\sinh^2\frac{(\alpha-\beta)\delta}2,\\
		&R_1=(p^2+\beta^2)^2-4\alpha\beta p^2,\quad R_2=(p^2+\beta^2)^2+4\alpha\beta p^2.
	\end{aligned}
	\label{22}
\eeq
The matrix $G(p,s)$ resembles its analogue in the steady-state case  (Antipov and Smirnov, 2013). However, in the steady-state 
problem, $\Ga$ and $\Gb$ are constants, while in the transient
case, they are functions of $p$ and $s$.
 
\setcounter{equation}{0}

\section{Transient problem for a plane}

\label{section:scalar}

In this section, we develop explicit representations for the solution and the weight functions of the transient problem on propagation of a semi-infinite crack in the 
particular case $\Gd=\infty$. Although the solution to this problem is known  (Freund, 1990), 
our solution has a different form. It is used  as a building block
for the approximate procedure proposed in the next section for
the solution of the problem on a crack
in a half-plane.

\subsection{Scalar RHP}
\label{section:scalar-problem}

On passing to the limit $\delta\to\infty$ in (\ref{22})
we arrive at the following two separate equations:
\beq
	\tilde\sigma^+_j(p,s)=\mu i g_j(p,s) \tilde\chi^-_j(p,s)+\tilde q_j^-(p,s),\quad p\in L,\quad j=1,2.
	\label{26}
\eeq
that are scalar RHPs with the coefficients
\beq
g_1(p,s)=\fr{R_1}{2\beta (p^2-\beta^2)p},
\quad g_2(p,s)=\fr{R_1}{2\Ga (p^2-\beta^2)p}.
\label{27}
\eeq
Assume first that $s$ is real and positive and let $s=c_ls'$ and $p=p's'$. Then we observe that
\beq
\Ga(p,s)= s'\tilde \Ga(p'), \quad \Gb(p,s)=s'\tilde \Gb(p'),
\label{28}
\eeq
where
$$
\tilde\Ga(p')=\hat\Ga(p'-a'_-)^{1/2}(p'-a'_+)^{1/2}, \quad
\tilde\Gb=\hat\Gb(p'-b'_-)^{1/2}(p'-b'_+)^{1/2},
$$
\beq
a'_\pm=\fr{i}{v_l\pm 1}\in{\mathbb C^\pm},
 \quad b'_\pm=\fr{i}{v_l\pm c_s/c_l}\in{\mathbb C^\pm}.
	\label{28.0}
\eeq
The branches of the functions $(p'-a'_\pm)^{1/2}$ and  $(p'-b'_\pm)^{1/2}$
are chosen as it is done in (\ref{13}) for the original functions. For simplicity, rename $p'$
as $p$  and write the 
problems (\ref{26}) as
\beq
	\tilde\sigma^+_j(ps',c_ls')=\mu i \tilde g_j(p) \tilde\chi^-_j(ps',c_ls')+\tilde q_j^-(ps',c_ls'),\quad p\in L,\quad j=1,2,
	\label{28.1}
\eeq
with the coefficients $\tilde g_j(p)$ independent of $s'$, 
\beq
\tilde g_j(p)=\fr{\tilde R_1}{2\tilde\beta (p^2-\tilde\beta^2)p},
\quad \tilde g_2(p)=\fr{\tilde R_1}{2\tilde\Ga (p^2-\tilde\beta^2)p},
\quad \tilde R_1=(p^2+\tilde\beta^2)^2-4\tilde\alpha\tilde\beta p^2.
\label{28.2}
\eeq
The functions $\tilde g_j(p)$ to be factorized have the following asymptotics at infinity and zero:
\beq
\tilde g_j(p)=\mp\Gg_j+O\left(\fr{1}{p}\right),\quad p\to\pm\infty, \quad  \tilde g_j(p)\sim -\fr{\Gg_j^\circ}{p}, \quad p\to 0,
\label{29}
\eeq
where
\beq
\Gg_1=\fr{R_0}{2\hat\Gb v_s^2}, \quad \Gg_2=\fr{R_0}{2\hat\Ga v_s^2}, \quad \Gg_1^\circ=\fr{1}{2c_s}, \quad \Gg_2^\circ=\fr{c_l}{2c^2_s},\quad
R_0=4\hat\alpha\hat\beta-(1+\hat\beta^2)^2.
\label{30}
\eeq
We emphasize that in the sub-Rayleigh regime, $V<c_R$, $c_R$ is the Rayleigh speed, and the parameter $R_0$ is positive.
Employing the relations (\ref{29}) we split the coefficients of the RHPs as
\beq
	\tilde g_j(p)=-\gamma_j\coth(\pi p) g_j^\circ(p).
	\label{31}
\eeq
Because the coefficients of the Riemann-Hilbert problem, $\tilde g_1$ and $\tilde g_2$, have a first order infinity at the point $p=0\in L$, to bypass this point, we deform the contour $L$.
There are two possibilities to do so. One of them is to replace $L$ by  $L_\Gve=L'\cup C^-_{\Gve}\cup L''$, where $L'=\{-\infty<p_1\le -\Gve, p_2=0\}$,
$L''=\{\Gve\le p_1< +\infty, p_2=0\}$, $C_\Gve^-=\{|p|=\Gve, p_2< 0\}$
, $p=p_1+ip_2$.
The second one is to replace $L$ by  $\tilde L_\Gve=L'\cup C^+_{\Gve}\cup L''$, where
 $C_\Gve^+=\{|p|=\Gve, p_2> 0\}$.
 Consider the first case. In Appendix A, we analyze the second possibility and show that the final solution is
 independent of the way we deform the contour.
The contour $L_\Gve$ splits the $p$-plane into two domains $D^+\ni 0 $ and $D^-$. It is important to notice that the new functions
 $ g_j^\circ(p)$ can be easily factorized
 $g_j^\circ(p)=\GO_j^+(p)/\GO_j^-(p)$, $p\in L_\Gve, $
  in terms of the Cauchy integrals
\beq
\GO_j^\pm(p)=\lim_{\tilde p\to p\in L_\Gve, \tilde p\in D^\pm}\GO_j(\tilde p), \quad	\Omega_j(\tilde p)=\exp\left\{\frac1{2\pi i}\int_{L_\Gve}\frac{\ln g_j^\circ(\tau)}{\tau-\tilde p}d\tau\right\},\quad \tilde p\in D^\pm,
	\label{32}
\eeq
since the functions $\ln g_j^\circ(\tau)$
are H\"older-continuous on the contour, $g^\circ_j(\tau)=1+O(1/\tau)$, $\tau\to\pm\infty$, positive at zero, and the increment of the argument
of $g_j^\circ(\tau)$ equals zero as $\tau$ traverses the whole contour $L_\Gve$.  

After factorizing the function $\coth(\pi p)$ in terms of the Gamma-functions
\beq
	\coth(\pi p)=\frac{iK^+(p)}{K^-(p)},\quad K^+(p)=\frac{\Gamma(1-i p)}{\Gamma(1/2-i p)},\quad K^-(p)=\frac{\Gamma(1/2+i p)}{\Gamma(i p)},
	\label{34}
\eeq
it is possible to transform the boundary condition (\ref{26}) of the RHP to the form
\beq
	\frac{\tilde\sigma^+_j(ps',c_ls')}{K^+(p)\Omega_j^+(p)}-\GY_j^+(p,s')=
	\frac{\mu\Gg_j\tilde\Gc^-_j(ps',c_ls')}{K^-(p)\Omega_j^-(p)}-\GY_j^-(p,s'), 
	\quad p\in L_\Gve,
	\label{35}
\eeq
where 
\beq
\GY_j^\pm(p,s')=\lim_{\tilde p\to p\in L_\Gve, \tilde p\in D^\pm}\GY_j(\tilde p,s'), \quad
\Psi_j(\tilde p,s')=\frac1{2\pi i}\int_{L_\Gve}\frac{\tilde q_j^-(\tau s',c_l s')}{K^+(\tau)\Omega_1^+(\tau)}\frac{d\tau}{\tau-\tilde p},\quad \tilde p\in D^\pm.
	\label{36}
\eeq
On applying the continuity principle and  the Liouville theorem and employing the asymptotics
$$
K^\pm(p)\sim (\mp ip)^{1/2}, \quad \GO_j^\pm(p)\sim 1,\quad \GY^\pm_j(p,s')=O(p^{-1}),
$$
\beq
\tilde\Gs_j^+(ps',c_l s')=O(p^{-1/2}), \quad \tilde\Gc_j^-(ps',c_ls')=O(p^{-1/2}), \quad p\to\infty,
\label{37}
\eeq
we find that the solution is unique and given by
\beq
	\begin{aligned}
		&\tilde\sigma^+_j(ps',c_ls')=K^+(p)\Omega_j^+(p)\Psi_j^+(p,s'),&& p\in D^+,\\
		&\tilde\chi^-_j(ps',c_ls')=(\mu\gamma_j)^{-1}K^-(p)\Omega_j^-(p)\Psi_j^-(p,s'),&& p\in D^-.
	\end{aligned}
	\label{38}
\eeq
Passage to the limit $\Gve\to0^+$ shows that the functions $\tilde\chi^-_j(-i\Gve s',c_ls')\to 0$ that is consistent with the fact that the displacement jumps  vanish as $x\to-\infty$.

\subsection{Stress intensity factors and the weight functions}
\label{section:scalar-sif}

We now turn our attention to the SIFs, $K_{I}(t)$ and $K_{II}(t)$,
determined as 
\beq
	\sigma_{xy}(x,0,t)\sim\frac{K_{II}(t)}{\sqrt{2\pi}}x^{-1/2},\quad \sigma_{yy}(x,0,t)\sim\frac{K_{I}(t)}{\sqrt{2\pi}}x^{-1/2},\quad x\to0^+.
\eeq
The same relations hold for the Laplace-transformed stresses and SIFs
\beq
	\hat\sigma_{xy}(x,s)\sim\frac{\hat K_{II}(s)}{\sqrt{2\pi}}x^{-1/2},\quad \hat\sigma_{yy}(x,s)
	\sim\frac{\hat K_{I}(s)}{\sqrt{2\pi}}x^{-1/2},\quad x\to 0^+.
\eeq
By the Abelian theorem  (Noble, 1988) and in view of (\ref{19}),
$$
	\tilde\sigma_1^+(ps',c_ls')\sim\tilde K_{II}(s')e^{i\pi/4}(2ps')^{-1/2},
	\quad 
	\tilde\sigma_2^+(ps',c_ls')\sim\tilde K_{I}(s')e^{i\pi/4}(2ps')^{-1/2},
	$$
	\beq
	 p\to\infty,\; \arg p\in(0,\pi).
	\label{40}
\eeq	
On the other hand, analysis of the solution (\ref{38}) shows that
\beq
	\tilde\sigma_j^+(ps',c_ls')\sim -e^{-i\pi/4}\GY_j^\circ(s')p^{-1/2},
\quad p\to\infty,\; \arg p\in(0,\pi),
	\label{41}
\eeq	
where
\beq
\GY_j^\circ(s')=\fr{1}{2\pi i}\int_L\fr{\tilde q_j^-(\tau s',c_ls')d\tau}{K^+(\tau)\GO_j^+(\tau)}.
\label{42}
\eeq
The consistency of formulas (\ref{40}) and  (\ref{41}) gives the desired expressions for the Laplace transforms of the SIFs
\beq
\hat K_{II}(s')=\sqrt{2} i\GY_1^\circ (s') (s')^{1/2},
\quad
\hat K_{I}(s')=\sqrt{2} i\GY_2^\circ (s') (s')^{1/2}.
\label{43}
\eeq
These expressions are obtained under the assumption that $s$ is real and positive.
By continuing analytically $\hat K_I(s)$ and $\hat K_{II}(s)$ from the real positive semi-axis
to the domain $\R s>0$ we define them in the whole right half-plane and may apply the inverse
Laplace transform to recover the SIFs
\beq
K_j(t)=\fr{1}{2\pi i}\int_\CL \hat K_j\left(\fr{s}{c_l}\right) e^{st}ds, \quad j=I,II,
\label{44}
\eeq
$\CL=\{\R s=\Gs>0, |\I s|<\infty\}$.

The SIFs can also be expressed through the weight functions,
 $W_{I}(\cdot,t)$ and $W_{II}(\cdot,t)$, by
\beq
	K_I(t)=\int_{-\infty}^{Vt} \sigma_{22}^\circ(x_1,0,t)W_I(x_1,t)dx_1,\quad K_{II}(t)=\int_{-\infty}^{Vt}\sigma_{12}^\circ(x_1,0,t)W_{II}(x_1,t)dx_1.
	\label{46}
\eeq
The weight functions $W_{I}(x_0,t)$ and $W_{II}(x_0,t)$ coincide with the SIFs determined for the special case 
of loading,
$\sigma_{j2}^\circ(x_1,0,t)=\delta(x_1-x_0)$, $j=1,2,$
where $\delta(x)$ is the Dirac $\delta$-function. In this case the functions $q_j(x,s)$ and
$\tilde q_j^-(p,s)$ introduced in (\ref{18}) and (\ref{19}) become
\beq
 q_j(x,s)=\fr{1}{V}e^{s(x-x_0)/V}, \quad \tilde q_j^-(ps',c_ls')=\fr{e^{-s' x_0/v_l}}{c_ls'
(1+i pv_l)}.
	\label{47}
\eeq
The integrals (\ref{42}) can be evaluated exactly
\beq
	\Psi_j^\circ(s') = \fr{e^{-sx_0/v_l}}{iVs'K^+(i/v_l)\GO_j^+(i/v_l)}, \quad j=1,2.
	\label{49}
\eeq
Consequently, the Laplace transforms of the weight functions have the form
\beq
	\hat W_I(x_0,s')=\frac{\sqrt{2}e^{-s' x_0/v_l}}{V\sqrt{s'}K^+(i/v_l)\Omega^+_2(i/v_l)},\quad 
	\hat W_{II}(x_0,s)=\frac{\sqrt{2}e^{-s' x_0/v_l}}{V\sqrt{s'}K^+(i/v_l)\Omega^+_1(i/v_l)},
	\label{50}
	\eeq
Now we can evaluate the inverse Laplace transforms exactly
\beq
W_{I}(x_0,t)=\sqrt{\fr{2}{\pi (Vt-x_0)}}w_{I},\quad
W_{II}(x_0,t)=\sqrt{\fr{2}{\pi (Vt-x_0)}}w_{II}
\label{50'}
\eeq
where
\beq
w_j=\fr{\GG(1/2+1/v_l)}{\sqrt{v_l} \GG(1+1/v_l)\GO_{3-j}^+(i/v_l)}, \quad j=I,II (1,2).
\label{50''}
\eeq
\begin{figure}
	\begin{center}\includegraphics{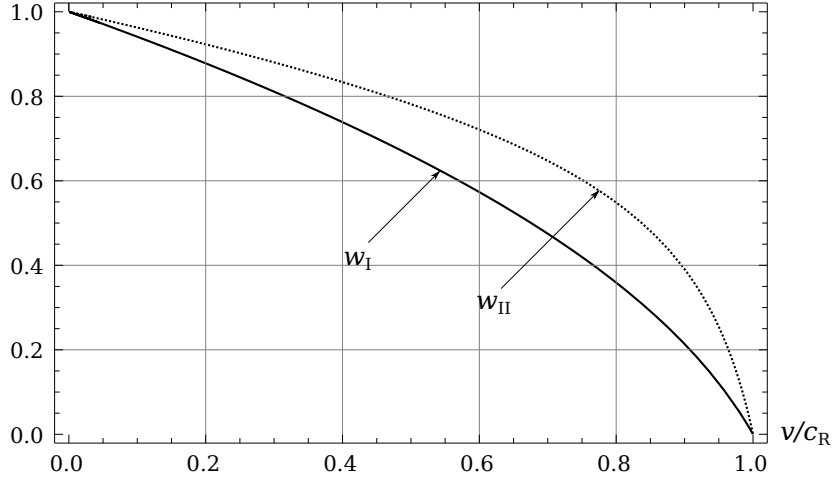}\end{center}
	\caption{Graphs of the functions $w_I$ and $w_{II}$ versus $V/c_R$ for $\nu=0.3$ in the case of an unbounded  plane.}
	\label{figure:plane}
\end{figure}
Graphs of the dimensionless functions $w_I$ and $w_{II}$
versus the dimensionless speed $V/c_R$ for $\nu=0.3$ is shown in Fig. \ref{figure:plane}. The graph
of the function $w_I$ is in good agreement with the one presented
by Freund (1990, p. 349).

\setcounter{equation}{0}

\section{Approximate solution of the transient problem for a half-plane}
\label{section:vector}

Since the structure of the matrix $G(p,s)$ given by (\ref{21}), (\ref{22})
does not allow for its explicit factorization by the methods currently available in the literature,
we propose a method of partial factorization. This technique
eventually leads to a system of two integral equations convenient 
for the determination of the SIFs and numerical implementation.

\subsection{System of integral equations}

We show first that the
direct use of the convolution theorem reduces the boundary condition of the vector RHP (\ref{20})
to a system of integral equations not convenient for numerical purposes. The system has the form 
\beq
	\int_{-\infty}^0 K(x-\xi,s)\left(\begin{array}{c}
	\hat\chi_1(\xi,s)\\\hat\chi_2(\xi,s)\end{array}\right)d\xi
	=-\frac1\mu\left(\begin{array}{c} q_1(x,s)\\ q_2(x,s)\end{array}\right),\quad -\infty<x<0,
	\label{57}
\eeq
where
\beq
	K(\xi,s)=\frac{i}{2\pi}\int_{-\infty}^{\infty}G(p,s) e^{-i p\xi}dp.
	\label{58}
\eeq
On analyzing 
the asymptotic behavior of the entries of $G(p,s)$
as $p\to 0$ and $p\to\infty$  we discover  that 
$$
g_{jj}(p,s)\sim -\Gg_j\sgn p\left[1+\fr{r_j}{p}+O\left(\fr{1}{p^2}\right)\right], \quad j=1,2,
$$
\beq
g_{12}(p,s)\sim r_0 e^{-2\hat\Gb \Gd|p|}, \quad p\to \pm\infty,
\label{61}
\eeq
where $r_j$ ($j=0,1,2$) are nonzero constants, and $\Gg_j$ are the positive
constants given by (\ref{30}). As $p\to 0$,
$g_{jj}(p,s)\sim -\tilde\Gg_j p^{-1}$,  $g_{12}(p,s)\sim -\tilde\Gg_0$,
($\tilde\Gg_j$ are positive constants). To clarify the structure of the diagonal kernels, 
we represent the functions $g_{jj}(p,s)$ as
\beq
g_{jj}(p,s)=-\Gg_j[\coth(\pi p)+g_{jj}^\circ(p,s)],
\label{64}
\eeq
where
\beq
g_{jj}^\circ(p,s)=\fr{r_j}{|p|}+O\left(\fr{1}{p^2}\right), \quad p\to\pm\infty,\quad g^\circ_{jj}(p,s)\sim\fr{\hat\Gg_j}{p}, \quad p\to 0,
\label{65}
\eeq
$\hat\Gg_j$ are constants. Because of the integral
\beq
\int_{-\infty}^\infty\coth(\pi p) e^{-i p\Gz}dp=-i\coth \fr{\Gz}{2},
\label{66}
\eeq
this ultimately brings us to
the system of singular integral equations
$$
\int_{-\infty}^0\left[\coth\fr{\Gx- x}{2}+k_{jj}(x-\Gx,s)\right]\hat\Gc_j(\Gx,s)d\Gx
$$
\beq
+\int_{-\infty}^0k_{j\,3-j}(x-\Gx,s)\hat\Gc_{3-j}(\Gx,s)d\Gx=-\fr{2\pi}{\mu\Gg_j}q_j(x,s), \quad -\infty<x<0, \quad j=1,2.
\label{70}
\eeq
The functions $k_{jj}(\Gz,s)$ have a logarithmic singularity, and the functions
$k_{j\, 3-j}(\Gz,x)$ are bounded at $\Gz=0$. As $\Gz\to\infty$, all the kernels decay, $k_{ij}=O(\Gz^{-1})$.
Difficulties will arise, however, if we try to use this system for computations. This is for the
slow convergence of the integrals in (\ref{70}) due to the presence
of the function $\coth\Gz/2$ bounded as $\Gz\to\infty$.

To avoid dealing with such kernels, we propose another approach.
First, we recast the system (\ref{57}) into a different form.
As in section \ref{section:scalar-problem}, we split the diagonal entries of the matrix $G$ 
\beq
	g_{jj}(p,s)=-\gamma_j\coth(\pi p) \check g_{jj}(p,s),\quad j=1,2,
	\label{71}
\eeq
factorize the function $\coth(\pi p)$ as in (\ref{34}) and the functions $\check g_{jj}(p,s)$ as follows:
$$
	\check g_{jj}(p,s)=\frac{\Omega_{jj}^+(p,s)}{\Omega_{jj}^-(p,s)},\quad p\in L_\Gve,
$$
	\beq
	 \Omega_{jj}(p,s)=\exp\left\{\frac1{2\pi i}\int_{L_\Gve}\frac{\ln \check g_{jj}(\tau,s)d\tau}{\tau-p}\right\}, \quad p\in D^\pm.
	\label{72}
\eeq
Introduce new functions
\beq
	\begin{aligned}
	&\check\sigma^+_j(p,s)=\frac{\tilde\sigma^+_j(p,s)}{K^+(p)\Omega_{jj}^+(p,s)},&\quad& \check\chi^-_j(p,s)=\frac{\mu\tilde\chi^-_j(p,s)}{K^-(p)\Omega_{jj}^-(p,s)},\\
	&\check q_j(p,s)=\frac{\tilde q_j^-(p,s)}{K^+(p)\Omega_{jj}^+(p,s)},&\quad& j=1,2.
	\end{aligned}
	\label{73}
\eeq
After the partial factorization has been implemented, in the new notations,
the original vector RHP (\ref{20}) reads
\beq
	\left(\begin{array}{c}\check\sigma_1^+(p,s)\\
\check\sigma_2^+(p,s)\end{array}\right)=\left(\begin{array}{cc}
\gamma_1&\check g_{1}(p,s)\\\check g_{2}(p,s)&\gamma_2\end{array}
	\right)\left(\begin{array}{c}\check\chi_1^-(p,s)\\
\check\chi_2^-(p,s)\end{array}\right)
+\left(\begin{array}{c}\check q_1(p,s)\\\check q_2(p,s)\end{array}\right),\quad p\in L_\Gve,
	\label{74}
\eeq
where
\beq
	\check g_{1}(p,s)=-\frac{ig_{12}(p,s)\Omega_{22}^-(p,s)}{\coth(\pi p)\Omega_{11}^+(p,s)},\quad
	 \check g_{2}(p,s)=\frac{ig_{12}(p,s)\Omega_{11}^-(p,s)}{\coth(\pi p)\Omega_{22}^+(p,s)},
\label{74'}
\eeq
and the functions $\GO_{jj}^\pm(p,s)$ are defined  by the Sokhotski-Plemelj formulas
\beq
\GO_{jj}^\pm(p,s)=\exp\left\{\pm\fr12\ln \check g_{jj}(p,s)+P.V.\GO_{jj}(p,s)\right\},\quad p\in L_\Gve.
\label{77'}
\eeq

Assume first that $x$ is negative. On applying the convolution theorem to (74)
we conclude that the vector RHP (\ref{74}) yields
$$
	\gamma_1\chi^*_1(x,s)+\int_{-\infty}^0 k_1^*(x-\xi,s)\chi_2^*(\xi,s)d\xi=- q_1^*(x,s),
	$$
	\beq
	\gamma_2\chi^*_2(x,s)+\int_{-\infty}^0 k_2^*(x-\xi,s)\chi_1^*(\xi,s)d\xi=- q_2^*(x,s),\quad -\infty<x<0.
\label{75}
\eeq
Here,
$$
\Gc_j^*(x,s)=\fr{1}{2\pi}\int_{L_\Gve}\check\Gc_j^-(p,s)e^{-ipx}dp,\quad
k_j^*(x,s)=\fr{1}{2\pi}\int_{L_\Gve}\check g_j (p,s)e^{-ipx}dp,
$$
\beq
q_j^*(x,s)=\fr{1}{2\pi}\int_{L_\Gve}\check q_j(p,s)e^{-ipx}dp,\quad j=1,2.
\label{76}
\eeq
Because of the asymptotics (\ref{61}) of the function $g_{12}(p,s)$,   the functions $\check g_j(p,s)$ decay exponentially as $p\to\pm\infty$.
It also follows from (\ref{74'}) that the functions $\check g_j(p,s)$ are continuously differentiable on the whole real axis and therefore  $|k_j^*(x,s)|\le cx^{-2}$
when  $x\to -\infty$ ($c$ is a function of $s$ and independent of $x$).
By the Riemann-Lebesgue lemma, the functions $\Gc_j^*(x,s)\to 0$ as $x\to-\infty$.

If  $x>0$, then the convolution theorem applied to (\ref{74})  gives 
$$
 	\int_{-\infty}^0 k_1^*(x-\xi,s)\chi_2^*(\xi,s)d\xi=\Gs_1^*(x,s)- q_1^*(x,s),
$$
\beq
	\int_{-\infty}^0 k_2^*(x-\xi,s)\chi_1^*(\xi,s)d\xi=\Gs_2^*(x,s)- q_2^*(x,s),\quad 0<x<\infty,
\label{75.2}
\eeq
where
\beq
\Gs_j^*(x,s)=\fr{1}{2\pi}\int_{L_\Gve}\check\Gs_j^+(p,s)e^{-ipx}dp,\quad j=1,2.
\label{75.3.0}
\eeq
On employing the continuity of the convolutions $k_1\ast\Gc_2^*$ and $k_2\ast\Gc_1^*$ 
and concatenating equations (\ref{75}) and (\ref{75.2}) at $x=0$ we establish the important relations
\beq
\Gg_j\Gc_j^*(0^-,s)=-\Gs_j^*(0^+,s)+q_j^*(0^+,s)-q_j^*(0^-,s), \quad j=1,2,
\label{75.3}
\eeq
to be used in the next section for computing the weight functions.

\subsection{Weight functions}

Because the mode-I and mode-II are coupled we have four weight functions, $W_{I,I}$, $W_{I,II}$, $W_{II,I}$, and $W_{II,II}$. Through them, the SIFs 
are found to be
$$
		K_I(t)=\int_{-\infty}^{Vt} W_{I,I}(x_1,t)\sigma_{22}^\circ(x_1,0)dx_1+\int_{-\infty}^{Vt}W_{I,II}(x_1,t)\sigma_{12}^\circ(x_1,0)dx_1,
		$$\beq
		K_{II}(t)=\int_{-\infty}^{Vt} W_{II,I}(x_1,t)\sigma_{22}^\circ(x_1,0)dx_1+\int_{-\infty}^{Vt}W_{II,II}(x_1,t)\sigma_{12}^\circ(x_1,0)dx_1.
	\label{75.4}
\eeq
The values of the weight functions $W_{I,I}$ and $W_{II,I}$ at a point $x_0$ coincide with the SIFs $K_I$ and $K_{II}$, respectively if $\sigma_{22}^\circ(x_1,0)=\delta(x_1-x_0)$
and $\sigma_{12}^\circ(x_1,0)=0$. Similarly, if $\sigma_{22}^\circ(x_1,0)=0$ and $\sigma_{12}^\circ(x_1,0)=\delta(x_1-x_0)$, then
the SIFs $K_I$, $K_{II}$ are equal to the other two weight functions $W_{I,II}(x_0,t)$ and $W_{II,II}(x_0,t)$, respectively.
As in section \ref{section:scalar-sif},  the transforms of the traction components, $\tilde\Gs_1^+(p,s)$ and  $\tilde\Gs_2^+(p,s)$, have the asymptotics
$$
	\tilde\sigma_1^+(p,s)\sim\tilde K_{II}(s)e^{i\pi/4}(2p)^{-1/2},
	\quad 
	\tilde\sigma_2^+(p,s)\sim\tilde K_{I}(s)e^{i\pi/4}(2p)^{-1/2},
	$$
	\beq
	 p\to\infty,\; \arg p\in(0,\pi).
	\label{77}
\eeq	
To express the asymptotics of these functions through the solution to the system of integral equations (\ref{75}),
we note that due to the relation (\ref{75.3}), the continuity of $q_j^*(x,s)$
at $x=0$ and formulas (\ref{75.3.0}), (\ref{76}) and (\ref{74'}) 
\beq
\check\Gs_j^+(p,s)\sim\fr{\Gg_j\Gc_j^*(0^-,s)}{ip}, \quad p\to\infty, \quad p\in D^+.
\label{78}
\eeq
Now we employ formulas (\ref{73}) to discover
\beq
\tilde\sigma_j^+(p,s)\sim -e^{i\pi/4}\Gg_j\Gc_j^*(0^-,s)p^{-1/2},
	\quad  p\to\infty, \quad p\in D^+.
\label{79}
\eeq
In view of (\ref{77}) this momentarily gives the SIFs
\beq\begin{aligned}
	\hat K_{I}(s)=-\sqrt{2}\gamma_2\chi^*_2(0^-,s),\quad \hat K_{II}(s)=-\sqrt{2}\gamma_1\chi^*_1(0^-,s).
	\label{92}
\end{aligned}\eeq
Therefore, to determine the Laplace transforms of the SIFs, we have to know the solution of the system (\ref{75}) at the point $x=0$ only. 
This system needs to be solved for the special right-hand side,
\beq
 q_1^*(x,s)=0,\quad q_2^*(x,s)=\frac{e^{s(x-x_0)/V}}{V K^+(i s/V)\Omega^+_{22}(i s/V,s)},
	\label{94}
\eeq
in the case of the weight functions $W_{I,I}$ and $W_{II,I}$ and for
\beq
 q_1^*(x,s)=\frac{e^{s(x-x_0)/V}}{V K^+(i s/V)\Omega^+_{11}(i s/V,s)},\quad q_2^*(x,s)=0,
	\label{94'}
\eeq
in the case  of the weight functions $W_{I,II}$ and $W_{II,II}$.

The SIFs are recovered from their Laplace transform by the inversion formula  (\ref{44}). The inversion
can be implemented
by computing one of the  real integrals
$$
K_j(t)=\fr{2e^{\Gs t}}{\pi}\int_0^\infty\R\{\hat K_j(\Gs+i\tau)\}\cos \tau t\,d\tau,
$$
\beq
K_j(t)=-\fr{2e^{\Gs t}}{\pi}\int_0^\infty\I\{\hat K_j(\Gs+i\tau)\}\sin\tau t\,d\tau,\quad j=I,II,
\label{94''}
\eeq
and the preference should be made to the one with the better rate of convergence.

\subsection{Numerical results}

Here we describe the numerical procedure for evaluation of the weight functions. We recall that the weight functions
coincide with the SIFs provided loading is chosen as it was described in section 4.2.
Due to formulas (\ref{92}) the Laplace transforms of the SIFs require the knowledge of the solution of the system of integral equations (\ref{75})
at the point $x=0$,  that is $\Gc_j^*(0^-,s)$, $j=1,2$. We find it convenient to map the system (\ref{75})  on the semi-infinite interval
into another one on the interval (-1,1). This is achieved by introducing the variables
\beq
	\xi=\frac{\xi'-1}{\xi'+1},\quad -1<\xi'<1,\quad x=\frac{x'-1}{x'+1},\quad -1<x'<1.
\label{99}
\eeq
The new system is easily seen to be
\beq
\Gg_j \CX_j(x',s)+\int_{-1}^1 \CK_{3-j}(x',\Gx',s)\CX_j(\Gx',s)d\Gx'=-\CQ_j(x',s), \quad -1<x'<1,\quad j=1,2,
\label{100}
\eeq
where
\beq
\CX_j(x',s)=\Gc_j^*(x,s), \quad \CK_j(x',\Gx',s)=\fr{2k^*_j(x-\Gx,s)}{(\Gx'+1)^2}, \quad \CQ_j(x',s)=q^*_j(x,s).
\label{100'}
\eeq
We note that due to the asymptotics of the original kernels   $k_j^*(x,s)=O(x^{-2})$ as  $x\to -\infty$, the new kernels $\CK_j(x',\Gx',s)$ are
bounded as $\Gx'\to -1$. 
This circumstance implies that the system (\ref{100}) can be approximately solved by using the collocation method with the collocation points $\Gx_k$ ($k=1,2,\dots,N$)
chosen to be
the zeros of the degree-$N$ Legendre polynomial $P_N(x)$.  The system of $2N$ linear algebraic equations associated with the system (\ref{100}) has the form
$$
\Gg_j \CX_j(x_n,s)+\sum_{k=1}^N v_k\CK_{3-j}(x_n,x_k,s)\CX_{3-j}(x_k,s)=-\CQ_j(x_n,s), 
$$\beq
 n=1,2,\ldots,N, \quad j=1,2,
\label{101}
\eeq
where $v_k$ are the Gauss-Legendre weights given by $v_k=2(1-x^2_k)^{-1}[P'_N(x_k)]^{-2}$.

The chief difficulty in the implementation of this procedure is the evaluation of the principal value of the  
integrals in (\ref{77'}), $P.V.\GO_{jj}(p,s)$. It is helpful to recast them as integrals over the arc $l=\{|p'|=1$, $\arg p'\in(-\pi/2.\pi/2)\}$
\beq
\GO_{jj}(p,s)=\exp\left\{\fr{1+p'}{2\pi i} P.V. \int_{l}\fr{\GG_j(\tau',s)d\tau'}{\tau'-p'}\right\},
\label{102}
\eeq
where
\beq
\GG_j(\tau',s)=\fr{\ln \check g_{jj}(\tau,s)}{1+\tau'}, \quad \tau'=\fr{1+i \tau}{1-i\tau}, \quad p'=\fr{1+ip}{1-ip}.
\label{103}
\eeq

Among numerous approximate formulas for the principal value of the Cauchy integral over a circle we choose the following one (Parton and Perlin, 1982, p.116):
\beq
\GO_{jj}(p,s)=\exp\left\{\fr{1+p'}{2M+1}\sum_{j=-M}^M\GG_j(e^{i\Gt_j},s)\left[
\fr12+\fr{i\sin\fr{M}{2}(\Gt-\Gt_j)\sin\fr{M+1}{2}(\Gt-\Gt_j)}
{\sin\fr12(\Gt-\Gt_j)}
\right]\right\},
\label{104}
\eeq
for being simple and proving a good accuracy. Here, $\Gt=-i\ln p'$, $\Gt_j=2\pi j/(2M+1)$.
\begin{figure}
	\begin{center}\includegraphics{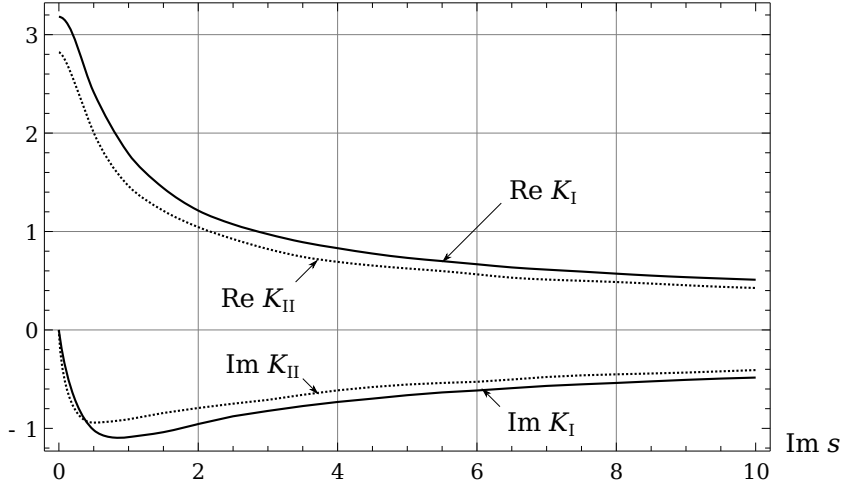}\end{center}
	\caption{Graph of the functions $\R K_I(s)$,  $\R K_{II}(s)$, $\I K_I(s)$, and $\I K_{II}(s)$
	for $\R s=0.5$.}
	\label{figure:hat}
\end{figure}

The final step in the evaluation of   the weight functions or, equivalently, the SIFs $K_I$ and $K_{II}$ with the special loads applied,
is the inversion of the Laplace transform. This can be done by applying one of the  formulas in (\ref{94''}). 
 For computations, we employ 
the uniform grid trapezoidal rule with $m+1$ grid points
 \beq
 K_j(t)\approx \frac{h e^{s_0 t}}{\pi}\left[\R\hat K_j(\Gs)+\R\hat K_j(\Gs+iT)\cos Tt+2\sum_{n=1}^{m-1}\R\hat K_j(\Gs+inh)\cos nht\right],
 \label{105}
 \eeq
 where $h$ is the grid spacing.
Our numerical results show (Fig. \ref{figure:hat}) that the rate of convergence is slow 
for both, the real and imaginary parts. To accelerate the convergence, we apply
the  Euler summation method for alternating series.
To transform (\ref{105}) into an alternating sum,  we put
 $h=\pi/(2t)$, $\Gs=A/(2t)$ and $T=\pi m/(2t)$, where $A$ is a fixed real positive constant. Then (Abate and Whitt, 1995)
\beq
	K_j(t)\approx \frac{e^{A/2}}{2t}\left[\R\hat K_j\left(\fr{A}{2t}\right)+\R\hat K_j\left(\fr{A+i\pi m}{2t}\right)\cos \fr{\pi m}{2}+2\sum_{n=1}^{m-1}(-1)^n\fr{\GD_n}{2^{n+1}}\right],
	\label{106}
	\eeq
	where
\beq	
	\Delta_n=\sum_{k=0}^n(-1)^k\left(\begin{array}{c}n\\k\end{array}\right)\R\left\{\hat K_j\left(\frac{A+2(n-k)\pi i}{2t}\right)\right\}.
\label{107}
\eeq
In our computations, following  Abate and Whitt (1995), we take $A=8\ln 10$.
\begin{figure}
	\begin{center}\includegraphics{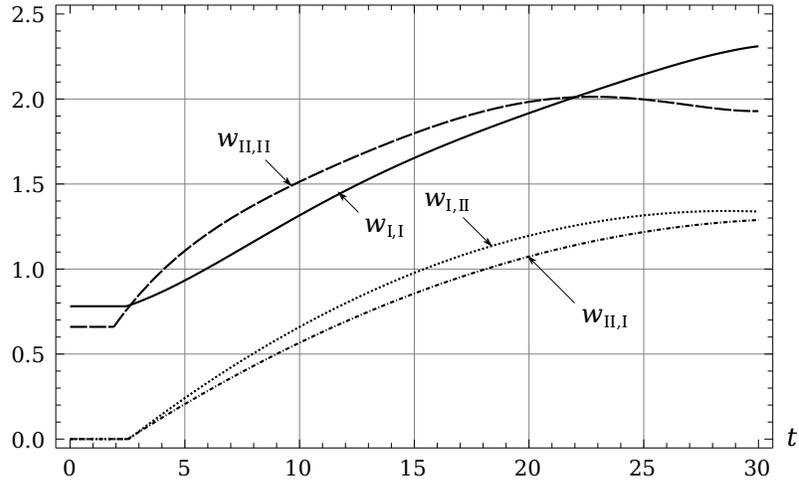}\end{center}
	\caption{The functions $w_{i,j}(0,t)=\sqrt{\fr12\pi V t}W_{i,j}(0,t)$  ($i,j=I,II$) versus time $t$ when	
	$\nu=0.3$, $\delta=1$\,m, $V=0.5c_R$\,m/s,  $c_l=1$\,m/s ($c_s\approx 0.5345$\,m/s, 
$c_R\approx0.4957$\,m/s).}
	\label{figure:KvsT}
\end{figure}

\begin{figure}
	\begin{center}\includegraphics{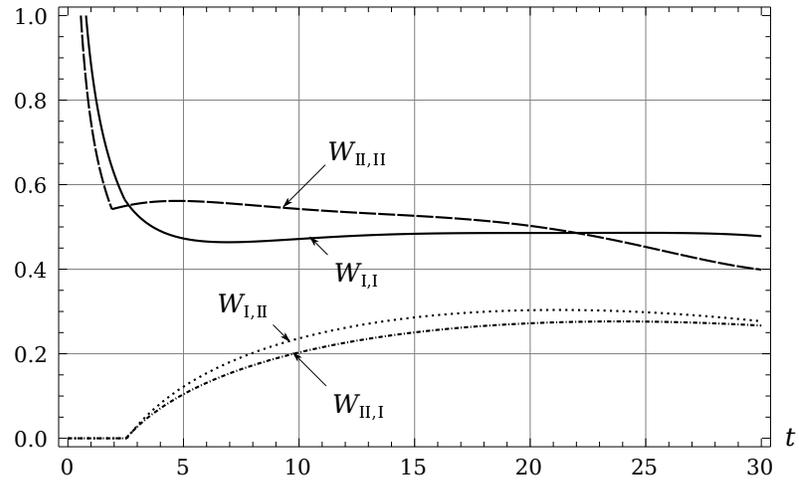}\end{center}
	\caption{
	The weight functions $W_{i,j}(0,t)$, ($i,j=I,II$) versus time $t$ when	
	$\nu=0.3$, $\delta=1$\,m, $V=0.5c_R$\,m/s,  $c_l=1$\,m/s ($c_s\approx 0.5345$\,m/s, 
$c_R\approx0.4957$\,m/s).}
\label{figure:KvsT2}
\end{figure}

\begin{figure}
	\begin{center}\includegraphics{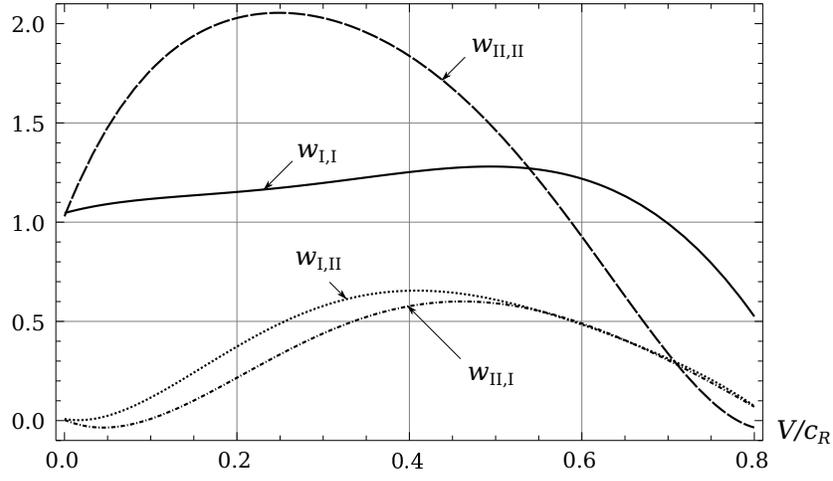}\end{center}
	\caption{The functions $w_{i,j}(0,t)$  ($i,j=I,II$)  versus  $V/c_R$ when		
	$\nu=0.3$, $\delta=1$\,m, $t=10s$,  $c_l=1$\,m/s ($c_s\approx 0.5345$\,m/s, 
$c_R\approx0.4957$\,m/s).}
\label{figure:KvsV}
\end{figure}

\begin{figure}
	\begin{center}\includegraphics{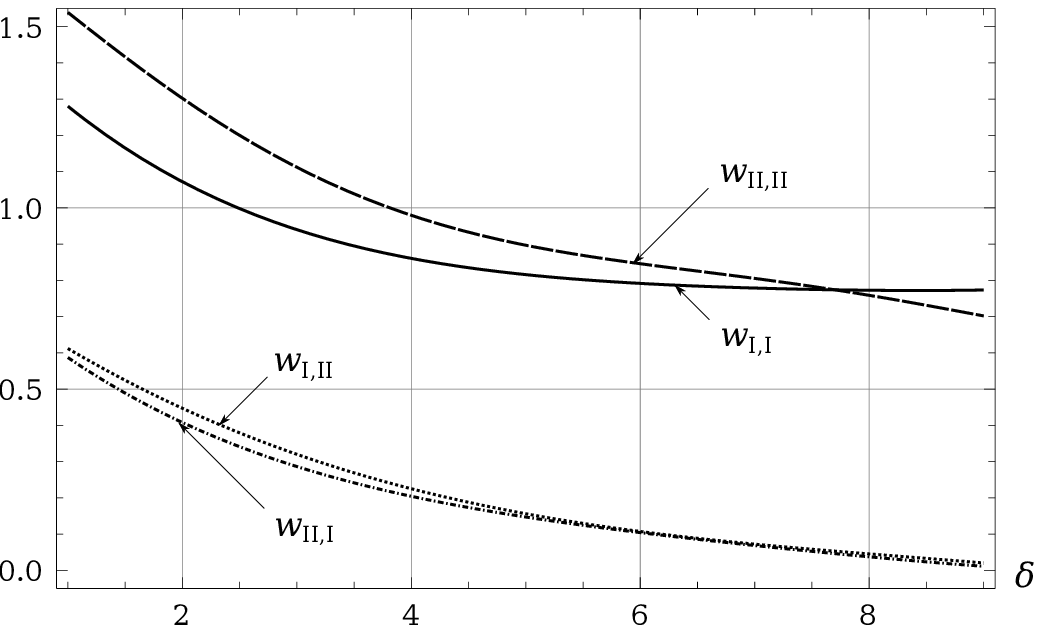}\end{center}
	\caption{
The functions $w_{i,j}(0,t)$  ($i,j=I,II$) versus the distance $\Gd$ from the
crack to the half-plane boundary 
when	$\nu=0.3$, $V=0.5c_R$\,m/s, $c_l=1$\,m/s, $t=10$\,s ($c_s\approx 0.5345$\,m/s, 
$c_R\approx0.4957$\,m/s).}
\label{figure:KvsD}
\end{figure}

Fig. \ref{figure:KvsT}  and \ref{figure:KvsT2} show how the functions $w_{i,j}(x_0,t)=\sqrt{\fr12\pi (V t-x_0)}W_{i,j}(x_0,t)$
and the weight functions $W_{i,j}(x_0,t)$ evolve with time. For computations, we  assume  $x_0$ to be zero, that is the time independent concentrated loads
are applied at the point $x_1=x_2=0$, the tip of the crack at time $t=0$. 
Since the material is stress-free for $t< 0$, it is expected that, when the crack starts propagating at constant speed $V$,
the elastic medium remains stress-free outside the disc of radius $c_l t$ centered at the  point $x_1= x_2=0$. At time $t_l'=\Gd/c_l$, the first longitudinal wave
strikes the boundary of the half-plane at the right angle, and at time $t=2t'_0$, it returns to 
the origin $x_1=y_1=0$. By that time, the crack tip has run the distance $2Vt_l'$,
and the distortion caused by the reflected wave reaches the crack tip
at time $t_l^*>2t_l'$ (for $\Gd>>1$, $t_l^*\sim 4t_1'$). The shear waves propagate slower, and the
corresponding time, when  the shear wave incident normally alters the SIFs,
is greater than $2\Gd/c_s>t_l^*$.
Due to other longitudinal waves reflected from the boundary at acute angles, 
the actual time when the boundary affects the SIFs may be smaller than $t_l^*$. The time when the
reflected longitudinal wave strikes  the crack at its tip
can be quickly evaluated. Let this wave   
hit the boundary of the half-plane  at time $t=t_l$ at angle $\Gt$ ($\Gt\in(\pi/2,\pi)$ is
measured from the incident wave direction to the boundary of the half-plane) (Fig. \ref{fig0}).
Then the reflected wave strikes the crack tip
at time $t=2t_l$. By that time, the crack has run the distance $2Vt_l$, and therefore,
$\sqrt{c_l^2t_l^2-\Gd^2}=Vt_l$. This implies 
\beq
t_l=\fr{\Gd}{\sqrt{c_l^2-V^2}}, \quad \Gt=\fr{\pi}{2}+\tan^{-1}\fr{1}{\sqrt{1/v_l^2-1}}.
\label{108}
\eeq
For the example used for drawing Fig. \ref{figure:KvsT}  and \ref{figure:KvsT2},
$\Gd=1\,$m, $V=0.5c_R$, and $c_R\approx 0.4957\,$m/s. Simple calculations show that
$2t_l\approx 2.0644\,$s and $\Gt=1.8213$. This time is consistent with the time $2t_l\approx 2\,$s
discovered from the approximate solution.
Our calculations (Fig. \ref{figure:KvsT}  and \ref{figure:KvsT2}) show that for time $0<t<2t_l$,
 the functions $w_{ii}(0,t)$  ($i=I,II$)  are constant and coincide with the parameters
$w_{i}$ associated with the mode-I and II weight functions for the whole plane and  given by (\ref{50''}).
The mixed mode functions $w_{I,II}(0,t)$ and $w_{II,I}(0,t)$
are identically equal to zero when $0<t<2t_l$.
The weight functions $W_{ij}(t,0)$ coincide with the corresponding weight functions of the problem on the whole plane for 
$0<t<2t_l$.   At time $t=2t_l$, the graphs of the weight functions associated with the half-plane and the plane start to diverge.

The functions  $w_{i,j}(0,t)$ versus the dimensionless speed $V/c_R$ are plotted in Fig. \ref{figure:KvsV}. As in the case of the whole plane,
the functions $w_{I,I}$ and $w_{II,II}$ tend to 1 and 0 when $V/c_R\to 0$ and $V/c_R\to 1$, respectively, while the off-diagonal functions, $w_{I,II}$
and $w_{II,I}$ tend to zero not only when $V/c_R\to 1$, but also when $V/c_R\to 0$. In the case of the whole plane,  the functions
$w_I$ and $w_{II}$ are monotonic, while in the case of the half-plane, they are not.

When the distance $\Gd$ from the crack to the boundary of the half-plane 
decreases, all the four functions $w_{i,j}(0,t)$ grow (see Fig. \ref{figure:KvsD}). As it is expected, when $\Gd\to\infty$, the functions
$w_{i,j}$ approach their limits, the corresponding functions for the whole plane, $w_{I,II}\to 0$,
$w_{II,I}\to 0$, and when $\nu=0.3$, $w_{I,I}\to w_I=0.781473$, $w_{II,II}\to w_{II}=0.659882$.

\setcounter{equation}{0}

\section{Crack growth at nonuniform speed beneath the boundary} 
\label{section:growth}

With the fundamental solution and weight functions at hand derived and computed in the previous
sections, we come now to the problem on nonuniform motion of a semi-infinite crack 
parallel to the boundary of a half-plane. We want to describe the motion of the crack
when speed, $V(t)$, is a prescribed smooth function for $t>0$. In order to do this,
we adjust the Freund (1990) approximate method proposed for a semi-infinite
crack moving at variable speed in an unbounded body. 

Suppose at time $t=0$ the crack starts moving, and its position at time $t$ is
described by $l(t)$, a continuously differentiable, nondecreasing function
such that $V(t)=l'(t)<c_R$. We approximate the curve $l(t)$ by a polygonal line
with the vertices ($t_k,l_k$), $l_k=l(t_k)$, $t_0=0$, $l_0=0$.
Denote $V_k=(l_{k+1}-l_k)/(t_{k+1}-t_k)$ the corresponding constant speed during the time $t_k<t<t_{k+1}$.

Initially, as $0<t<t_1$, the crack extends at speed $V_0=\const$ by negating the stresses $\Gs_{12}^0(x_1,0)$
and $\Gs_{22}^0(x_1,0)$ for $x_1>0$. They are determined from the solution of the static
problem, $P_{-1}$, on a semi-infinite crack parallel to the boundary of a half-plane.
This problem provides the starting point for a complete description of the 
nonuniform motion of the crack. 
An exact method of matrix Wiener-Hopf factorization for this problem was presented
by Zlatin and Khrapkov (1986) in the case when the forces were applied to 
the strip at infinity, and the boundary was free of traction. These authors  reduced the problem
to a homogeneous order-2 vector RHP, solved it exactly and found
the SIFs. On employing their method
it is possible to derive the exact solution of the inhomogeneous 
RHP for general loading and determine the stresses everywhere in the body including
the line $x_2=0$ ahead of the crack. We accept that the solution to Problem
$P_{-1}$ is already available.

Coming back now to the problem on a moving crack assume that the crack suddenly stops at time $t=t_1$ at 
the point $x_1=l_1, x_2=0$. Obviously, some stresses, $\tilde\Gs_{12}^1(x_1,0)$ 
 and $\tilde\Gs_{22}^1(x_1,0)$,
 are radiated out along the line $x_2=0$, $x_1>l_1$. 
These stresses are unknown
 {\it a priori}
and must be determined. To continue its motion, the crack negates these unknown stresses. This results in vanishing the SIFs when $x_1=V_0t>l_1$,
\beq
K_I(t;V_0)=0, \quad K_{II}(t;V_0)=0, \quad V_0t>l_1,
\label{5.1}
\eeq
and a necessity of solving a transient problem, $P_0$, arises. The statement of Problem  $P_0$
coincides with 
that given in Section 2 with the exception that $V=V_0$ and the boundary conditions (\ref{1.0})
on the faces of the crack read 
\beq
\Gs_{j2}=-\Gs_{j2}^0(x_1,0)\Gc_{(0,l_1)}(x_1)+\tilde\Gs_{j2}^1(x_1,0)\Gc_{(l_1, V_0t)}(x_1),
\quad -\infty<x_1<V_0t, \quad x_2=0^\pm,
\label{5.2}
\eeq 
 with $\tilde\Gs_{j2}^1(x_1,0)$ to be recovered from equations (\ref{5.1}). Here,
$\Gc_{(a_1,a_2)}(x_1)=1$ if $ x_1\in(a_1,a_2)$ and
$\Gc_{(a_1,a_2)}(x_1)=0$ otherwise.
To solve  equations (\ref{5.1}), we note the following remarkable property of the weight functions  
\beq
W_{i,j}(x_0,t;V)=W_{i,j}(0,t-x_0/V;V), \quad i,j=I,II.
\label{5.3} 
\eeq
To show this, we recall that due to (\ref{47}) the Laplace transforms of the loads for the weight functions are given by 
\beq
q_j(x,s)=e^{-sx_0/V}\fr{e^{sx/V}}{V}, \quad j=1,2,
\label{5.4}
\eeq
and  consequently from (\ref{19}), (\ref{73}) and  (\ref{76})  we derive
the relations
$$
\tilde q^-_j(p,s;x_0)=e^{-sx_0/V}\tilde q^-_j(p,s;0),
\quad
\check q_j(p,s;x_0)=e^{-sx_0/V}\check q^-_j(p,s;0),
$$
\beq
\check \Gc_j^-(p,s;x_0)=e^{-sx_0/V}\check \Gc^-_j(p,s;0),\quad
\Gc_j^*(x,s;x_0)=e^{-sx_0/V}\Gc^*_j(x,s;0),
\label{5.5}
\eeq
The latter formula and also (\ref{92}) imply that the Laplace transforms of the
weight functions satisfy the equation
\beq
\hat W_{i,j}(x_0,s;V)=e^{-sx_0/V}\hat W_{i,j}(0,s;V), \quad i,j=I,II.
\label{5.6}
\eeq
and the relation (5.3) holds. Combining these results we can write down formulas (\ref{75.4}) for 
the SIFs in the form
$$
		K_I(t;V_0)=-K'(t;V_0)	
		+\int_{l_1}^{V_0t} [W_{I,I}(0,t-x_1/V_0;V_0)
\tilde\Gs_{22}^1(x_1,0)	
$$$$
+W_{I,II}(0,t-x_1/V_0;V_0)\tilde\Gs_{12}^1(x_1,0)]dx_1,
		$$
		$$
K_{II}(t;V_0)=-K''(t;V_0)+\int_{l_1}^{V_0t} [W_{II,I}(0,t-x_1/V_0;V_0)
\tilde\Gs_{22}^1(x_1,0)	
$$
\beq
+W_{II,II}(0,t-x_1/V_0;V_0)\tilde\Gs_{12}^1(x_1,0)]dx_1,
\label{5.7}
\eeq
where $K'(t;V_0)$ and $K''(t;V_0)$ are known functions
$$
K'(t;V_0)=\int_{0}^{l_1} [W_{I,I}(0,t-x_1/V_0;V_0)
\Gs_{22}^0(x_1,0)		
+W_{I,II}(0,t-x_1/V_0;V_0)\Gs_{12}^0(x_1,0)]dx_1,				
$$
\beq
K''(t;V_0)=\int_{0}^{l_1} [W_{II,I}(0,t-x_1/V_0;V_0)
\Gs_{22}^0(x_1,0)		
+W_{II,II}(0,t-x_1/V_0;V_0)\Gs_{12}^0(x_1,0)]dx_1.
\label{5.8}
\eeq		
 We see now that the property (\ref{5.3}) allows for an exact solution of the system (\ref{5.1}) by transforming it into a system of two Volterra convolution equations
and applying the Laplace transform. Indeed, with a change of the variables
\beq
x_1=V_0\tau'+l_1, \quad t=\tau+l_1/V_0,
\label{5.9}
\eeq
and the functions to be found
\beq
\pi_I(\tau')=\tilde\Gs_{22}^1(V_0\tau'+l_1,0), \quad 
\pi_{II}(\tau')=\tilde\Gs_{12}^1(V_0\tau'+l_1,0),
\label{5.10}
\eeq
the system (\ref{5.1}) reads
\beq
\sum_{j=I}^{II}\int_0^\tau W_{i,j}(0,\tau-\tau';V_0)\pi_i(\tau')d\tau'=\Go_i(\tau), \quad \tau>0,\quad i=I,II,
\label{5.11}
\eeq
where
\beq
\Go_I(\tau)=V_0^{-1}K'(\tau+l_1/V_0;V_0), \quad  \Go_{II}(\tau)=V_0^{-1}K''(\tau+l_1/V_0;V_0).
\label{5.12}
\eeq
The Laplace images $\hat\pi_i(s)$ of the unknown functions $\pi_i(\tau')$
can be easily recovered from the system of linear algebraic equations
\beq
\sum_{j=I}^{II}\hat W_{i,j}(0,s;V_0)\hat\pi_i(s)=\hat \Go_i(s), \quad i=I,II.
\label{5.13}
\eeq
On performing the Laplace inversion we obtain
$$
\pi_I(\tau')=\fr{1}{2\pi i}\int_{\Gs-i\infty}^{\Gs+i\infty}
\fr{\hat W_{II,II}(0,s;V_0)\hat\Go_I(s)-\hat W_{I.II}(0,s;V_0)\hat\Go_{II}(s)}{\hat W(s;V_0)}e^{s\tau'}ds,
$$
\beq
\pi_{II}(\tau')=\fr{1}{2\pi i}\int_{\Gs-i\infty}^{\Gs+i\infty}
\fr{\hat W_{I,I}(0,s;V_0)\hat\Go_{II}(s)-\hat W_{II.I}(0,s;V_0)\hat\Go_{I}(s)}{\hat W(s;V_0)}
e^{s\tau'}ds,
\label{5.14}
\eeq
where $\Gs>0$ and 
\beq
\hat W(s;V_0)=\hat W_{I,I}(0,s;V_0)\hat W_{II,II}(0,s;V_0)-\hat W_{I,II}(0,s;V_0)\hat W_{II,I}(0,s;V_0).
\label{5.15}
\eeq
Thus the stresses to be nullified for $x_1>l_1$ have the form
\beq
\tilde\Gs_{22}^1(x_1,0)=\pi_I\left(\fr{x_1-l_1}{V_0}\right), \quad 
\tilde\Gs_{12}^1(x_1,0)=\pi_{II}\left(\fr{x_1-l_1}{V_0}\right),\quad x_1>l_1.
\label{5.15'}
\eeq 
We note that the Laplace transforms $\hat W_{i,j}(0,s;V_0)$ have already been determined.
They are expressed through the solution
at the point $0$ of the system of integral equations (\ref{75.2}) by (\ref{92})
as $V=V_0$
with the loading $\Gs_{22}^\circ(x_1,0)=\Gd(x_1)$ and $\Gs_{12}^\circ(x_1,0)=0$
for the weight functions $W_{I,I}(0,t;V_0)$ and  $W_{II,I}(0,t;V_0)$ and with
$\Gs_{22}^\circ(x_1,0)=0$ and $\Gs_{12}^\circ(x_1,0)=\Gd(x_1)$
for the functions $W_{I,II}(0,t;V_0)$ and  $W_{II,II}(0,t;V_0)$. 

In addition to nullifying the stresses $\tilde\Gs^1_{j2}(x_1,0)$, $j=1,2$, the solution for a suddenly stopped crack has to generate zero displacement jumps
through the line $x_2=0$
on the segment $l_1<x_1<V_0t$. In contrast to the whole plane problem, when  this
is possible to verify analytically  (Freund (1990)  for the sub-Rayleigh speeds and Huang and Gao (2002)
for the transonic regime), it is not visible how it can be done without deploying computer based computations.
That is why this condition needs to be tested numerically when the algorithm is applied.

By employing this procedure for the next period of time, $t_1<t<t_2$,
and determining the weight functions associated with speed $V=V_1$
we can find the loads  $\tilde\Gs_{i2}^2(x_1,0)$ ($i=1,2$)
needed to negate the stresses generated by the crack when it suddenly stops
at the point $x_1=l_2$. The boundary conditions (\ref{5.2}) for the corresponding 
problem $P_1$ read
\beq
\Gs_{j2}=-\Gs_{j2}^1(x_1,0)\Gc_{(l_1,l_2)}(x_1)+\tilde\Gs_{j2}^2(x_1,0)\Gc_{(l_2, V_1t)}(x_1),
\quad -\infty<x_1<V_1t, \quad x_2=0^\pm,
\label{5.16}
\eeq 
where the traction components $\Gs_{j2}^1(x_1,0)$ are known 
\beq
\Gs_{j2}^1(x_1,0)=
\Gs_{j2}^0(x_1,0)+\tilde\Gs_{j2}^1(x_1,0),
\label{5.17}
\eeq
while the components $\tilde\Gs_{j2}^2(x_1,0)$
have to be recovered from the system of two equations
$K_I(t;V_1)=0$, $ K_{II}(t;V_1)=0$, $V_1t>l_2$, that is equivalent to 
the corresponding system of two Volterra equations solvable by the Laplace transform 
as in the previous step.

Following the pattern established above, this procedure can be continued further up to any period of time $(t_{k},t_{k+1})$. It 
gives an approximate solution of the problem on motion of a semi-infinite
crack beneath the boundary at piecewise constant speed $V=V_i$, $t\in(t_{i},t_{i}+1)$, $i=0,1,\ldots,k$,
that approximates the original smooth function $V(t)$. The solution of this
model problem, $P$, is obtained 
by summing up the solutions of all Problems $P_i$ ($i=-2,-1,0,1,\ldots,k$), where
$P_{-2}$ is the elementary problem on a half-plane without a crack with somehow
prescribed traction on the boundary and internally loaded; its exact solution is 
available in the literature, $P_{-1}$ is the static problem for the semi-infinite crack with 
the traction components on the crack faces being prescribed such that they negate
the corresponding stresses coming from Problem $P_{-2}$. Problems $P_i$ ($i=0,1,\ldots,k-1$)
are the transient problems with the boundary conditions chosen accordingly.  The last problem $P_{k}$ is also a transient  problem whose boundary conditions pattern  is different from $P_i$ ($i=0,1,\ldots,k-1$) since they do not have 
stresses to be determined from the solution on a suddenly stopped crack.  The boundary conditions have the form
 \beq
\Gs_{j2}=-\left(\Gs_{j2}^0(x_1,0)+\sum_{i=1}^{k}\tilde\Gs_{j2}^i(x_1,0)\right)\Gc_{(l_{k},V_{k}t)}(x_1),
\quad -\infty<x_1<V_{k}t, \quad x_2=0^\pm,
\label{5.18}
\eeq 
Clearly, for the total problem $P$, the homogeneous boundary conditions
on the crack faces $\{0<x_1<l_{k}, x_2=0^\pm\}$ are satisfied. As for the SIFs at the tip of the crack at time 
$t\in(t_{k},t_{k+1})$, when the crack moves at speed $V_{k}$,
in general, they do not vanish, depend on time  and are defined by the SIFs generated by Problem 
$P_{k}$. 

A feature of Problem $P$ is in the presence of the boundary. As it was pointed out in the previous section, initially, when $t<2t_l$ ($t_l$ is given by (\ref{108})), and when
the longitudinal wave reflected from the half-plane boundary has not reached the crack, the off-diagonal weight functions
$W_{I,II}$ and $W_{II,I}$ vanish, and the diagonal functions $W_{I,I}$ and $W_{II,II}$
coincide with those associated with the problem on the whole plane with a crack.
Therefore, for this short period of time, the algorithm proposed
by Freund (1990) can be repeated without any changes. However, this does not mean that
the actual motion of the crack in a half-plane will be the same as in the whole plane
even  
for time $t<2t_l$. To make this conclusion, we need to recall that 
the boundary conditions of Problem $P_0$ depend on the stresses $\Gs_{i2}^0(x_1,0)$ ($i=1,2$) generated by the static crack in the half-plane which are apparently not the 
same as the ones associated with the whole plane.
When time exceeds $2t_l$, then, in general, all the weight functions $W_{i,j}$ are 
nonzero and different
from those associated with the whole plane. In this case to describe
the nonuniform crack motion, the algorithm
we exposed needs to be applied.

\section{Conclusion} 

We have derived the fundamental solution and the weight functions of the transient two-dimensional problem
on a semi-infinite crack propagating at constant speed parallel to the boundary of a half-plane.
The boundary of the half-plane is free of traction, while the crack faces are subjected to general  time-independent loading. 
We have reduced the boundary-value problem to a vector RHP on the real axis. In the particular case, when the crack
is far away from the boundary of the half-plane, the vector RHP is decoupled and solved by quadratures.
In the general case, we have split the matrix coefficient into a discontinuous diagonal matrix and a continuous matrix, factorized
the discontinuous part and rewritten the vector RHP as a system of two convolution equations on the segment $-\infty<x<0$. For numerical purposes,  it was recast as a system of two Fredholm integral equations on the  segment $(-1,1)$. 
We have derived the Laplace transforms of the SIFs and the weight functions in terms of the solution of the convolution equations
at the point $x=0$.  The Laplace transform has been inverted numerically. To improve the convergence, we  have
applied the Euler summation method for  alternating series. We have obtained numerical results for the SIFs when concentrated
loads are applied to the crack faces (at time $t=0$ at the crack tip).  This model problem generates four weight functions $W_{i,j}$, $i,j=I,II$. It has been discovered that
during a certain initial period of time, $0<t<2t_l$, the off-diagonal weight functions $W_{i,j}$, $i\ne j$,
vanish, and the diagonal functions coincide with the ones for the case of the whole plane.  For time $t>2t_l$, the boundary effects  play a significant role, and, in general, all the four weight functions
do not vanish and are different from  the corresponding functions associated with the whole plane.  
It has also been found that the dimensionless functions $w_{i,i}(0,t)=\sqrt{\fr12\pi V t}W_{i,i}(0,t)$ ($i=I,II$) tend to 1 and 0
as $V/c_R$ tends to 0 and 1, respectively ($V$ is the crack speed and $c_R$ is the Rayleigh speed), while $w_{i,j}$ ($i\ne j$)
vanish when $V/c_R$ approach both points, $0$ and 1.  We have found that $w_{ij}$ are not monotonic functions of $V/c_R$
and attain their local maximum in the interval $(0,V/c_R)$.
As the distance $\Gd$ from the crack to the boundary
decreases, all the functions $w_{ij}$ grow. We emphasize that apart from small $\Gd$ our numerical method is stable for 
all parameters $\Gd$.

Based on the approximate algorithm by Freund (1990) for the problem on a semi-infinite
crack propagated at a nonuniform rate in the whole plane
we have developed a procedure for the case when the crack propagates also at  prescribed variable sub-Rayleigh
speed in a half-plane parallel to the boundary and when the boundary effects
are significant. 
The implementation of the method requires solving a system of Volterra convolution equations
whose kernels are the associated weight functions, not 
a single Abel integral equation as in the whole plane case.
The  system of Volterra equations also admits a closed-form solution. However,
in the case of a half-plane,
there is no analog of the remarkable formula for the Mode I SIF  
$K_I(l(t), V_k)=k(V_k) K_I(l(t),0)$ in any interval $t_k<t<t_{k+1}$ derived  for the 
whole plane (Freund, 1990).
There is another difference between the whole plane and half-plane solutions. The displacement 
jumps though the crack line $x_2=0$
have to vanish on the segments $l_i<x_1<V_{i-1}t$, $i=1,\ldots, k$. This property was analytically
proved by Freund (1990) in the sub-Rayleigh regime and by Huang and Gao (2002) 
in the transonic regime. For the half-plane problem, this condition needs to be verified numerically
for each Problem $P_i$ ($i=0,1,\ldots,k-1$)
during the implementation of the procedure.

To compute the SIFs at time $t$, $2t_l<t_k<t<t_{k+1}$, for the crack in a half-plane,
one needs to derive the weight functions for all intermediate speeds $V_i$. 
We have shown that initially, before the longitudinal wave
reflected from the boundary strikes the crack and when the weight functions
coincide with those for the whole plane, the relatively simple Freund's
algorithm works. At the same time, the solution is still different since it relies
on the static solution on a cracked half-plane, not the whole plane with the crack.
When the first longitudinal wave reflected from the half-plane boundary reaches the crack surface moving at speed
$V(t)<c_R$, the boundary substantially affects the weight functions.
To determine the SIFs at the crack tip at some time $t\in (t_k,t_{k+1})$,
consequently, one may employ the procedure presented that requires solving
 the same transient problem for different 
constant speeds $V_i$ ($i=0,1,\ldots,k$)
and a system of Volterra equations to determine at each step the 
loads need to be negated to make possible for the crack to advance.

\appendix

\setcounter{equation}{0}

\section{Solution for a plane in the case of the contour $\tilde L_\Gve$}

It is important to show that as $\Gve\to 0^+$
the weight functions are invariant of the way the original contour $L$ is deformed. 
The contour $\tilde L_\Gve$  splits the $p$-plane into the domains $\tilde D^-\ni 0$ and $\tilde D^+$. 
In this case, $\coth(\pi p)$ in the representation (\ref{31}) needs to be factorized as follows:
\beq
	\coth(\pi p)=\frac{i\tilde K^+(p)}{\tilde K^-(p)},\quad \tilde K^+(p)=-\frac{\Gamma(-i p)}{\Gamma(1/2-i p)},\quad \tilde K^-(p)=\frac{\Gamma(1/2+i p)}{\GG(1+i p)}.
\label{51}
\eeq
Due to the fact that the asymptotics of the factors $\tilde K^\pm(p)$ at infinity is different
from that of $K^\pm(p)$, $\tilde K^\pm(p)\sim\mp(\mp p)^{-1/2}$, $p\to\infty$, $p\in \tilde D^\pm$,
the solution to the RHPs has arbitrary constants $C_j$,
\beq\begin{aligned}
	\tilde\sigma^+_j(ps',c_ls')&=\tilde K^+(p)\Omega_j^+(p)\left[C_j+\tilde\Psi_j^+(p,s')\right],\quad p\in\tilde D^+,\\
	\tilde\chi_j(ps',c_ls')&=(\mu\gamma_j)^{-1}\tilde K^-(p)\Omega^-_j(p) \left[C_j+\tilde\Psi_j^-(p,s')\right],\quad p\in\tilde D^-,\\
	\tilde\Psi_j^\pm(p,s')&=\frac1{2\pi i}\int_{\tilde L_\Gve}\frac{\tilde q^-_j(\tau s',c_ls')}{\tilde K^+(\tau)\Omega_j^+(\tau)}\frac{d\tau}{\tau-p},\quad p\in\tilde D^\pm,\quad j=1,2.
	\label{53}
\end{aligned}\eeq
Now, let $p=i\Gve\in\tilde L_\Gve$ ($\Gve>0$). 
According to the Sokhotski-Plemelj formulas
\beq
	\tilde\Psi^+_j(i\Gve,s')-\tilde\Psi_j^-(i\Gve,s')=\frac{\tilde q_j^-(i\Gve s',c_ls')}{\tilde K^+(i\Gve)\Omega_j^+(i\Gve)}.
\label{54}
\eeq
As $\Gve\to 0^+$, $[\tilde K^+(i\Gve)]^{-1}\to 0$, while the other two functions in the right-hand
side of (\ref{54}) are bounded and nonzero. Therefore, $\tilde\GY^\pm_j(i\Gve,s')\sim
\tilde\GY_j(0,s')$,  $\Gve\to 0^+$, where $\tilde\GY_j(0,s')$ is the principal value of the Cauchy integral
\beq
\tilde\Psi_j(p,s')=\fr{1}{2\pi i}
\int_{L}\frac{\tilde q^-_j(\tau s',c_ls')}{\tilde K^+(\tau)\Omega_j^+(\tau)}\frac{d\tau}{\tau-p},\quad p\in L,
\quad j=1,2.
\label{55}
\eeq 
It is evident that $\tilde\Gs_j^+(ps',c_ls')\to\infty$ as $p=i\Gve\to 0$ unless $C_j=-\tilde\GY_j(0,c_ls')$.
On the other hand, this choice of the constants $C_j$ guarantees that $\tilde\Gc_j^-(0,c_ls')=0$, and
both displacement jumps vanish at $x=-\infty$. Simple calculations show that
$\tilde\GY_j(0,s')=i\tilde\GY_j^\circ(s')$, where $\GY_j^\circ(s')$ is given by 
(\ref{42}). Analysis of $\tilde\Gs_j^+(ps',c_ls')$ as $p\to\infty$, $p\in\tilde D^+$, results in the asymptotics
\beq
\tilde\Gs_j^+(ps',c_ls')\sim(-ip)^{-1/2}C_j=-e^{-i\pi/4}\GY_j^\circ(s')p^{-1/2}, 
\quad p\to\infty, \quad \arg p\in(0,\pi),
\label{56}
\eeq
that coincides with formula (\ref{41}), and brings us  to the  expressions for the SIFs and the weight 
functions derived in section \ref{section:scalar-sif} in the case of the contour $L_\Gve$.

\vspace{.2in}

{\large {\bf References}}

\vspace{.1in}

Abate, J.,  Whitt, W.,  1995. Numerical inversion of Laplace transforms of probability distributions, ORSA J. Comput. 7, 36-43.

Antipov, Y.A., 1999. An exact solution of the 3D-problem of an interface semi-infinite plane crack,
J. Mech. Phys. Solids 47,1051-1093.

Antipov, Y.A. 2012. Weight functions of a crack in a two-dimensional micropolar solid, Quart. J. Mech. Appl. Math. 65, 239-271.

Antipov, Y.A., Moiseyev, N.G. 1991. Exact solution of the plane problem for a composite plane
with a cut across the boundary between two media, J. Appl. Math. Mech. (PMM) 55,  531-539.

Antipov, Y.A., Silvestrov, V.V. 2002.
Factorization on a Riemann surface in scattering theory, 
Quart. J. Mech. Appl. Math. 55 (4), 607-654.

Antipov, Y.A.,   Smirnov, A.V., 2013. 
Subsonic propagation of a crack parallel to the boundary of a half-plane, Math. Mech. Solids 18, 153-167.

Antipov, Y.A.,  Willis, J.R. 2003. Transient loading of a rapidly-advancing Mode-II crack in a viscoelastic medium, Mechanics of Materials 35, 415-431.

Antipov, Y.A.,  Willis, J.R. 2007. Propagation of a Mode-II crack in a viscoelastic medium with different bulk and shear relaxation, J. Engrg. Math. 59, 359-371.

Bercial-Velez, J.P., Antipov, Y.A.,   Movchan, A.B., 2005.  High order asymptotics and perturbation problems for 3D interfacial cracks, J. Mech. Phys. Solids 53, 1128-1162.

Bojarski, B.V. 1958 On the stability of the Hilbert problem 
for an unknown vector, Soob\u{s}\u{c}. Akad. Nauk Gruzin. SSR 25, 391-398.

Bueckner, H.F., 1970. A novel principle for the computation of stress intensity factors, Zeit. Angew.
Math. Mech.  50, 529-546.

Freund,  L.B., 1990.  Dynamic Fracture Mechanics, Cambridge University Press, Cambridge.

Gohberg, I. C., Kre\u in, M.G. 1958. On the stability of a system of partial indices of the 
Hilbert problem
for several unknown functions, Dokl. AN SSSR 119,  854-857.

Huang, Y.,  Gao, H., 2001.  Intersonic crack propagation - Part I: The fundamental solution, J. Appl. Mech. 68, 169-175.

Huang, Y.,  Gao, H., 2002. Intersonic crack propagation - Part II: Suddenly stopping crack, J. Appl. Mech. 69, 76-80.

Litvinchuk, G.S.,  Spitkovski\u i, I.M. 1987. Factorization of Measurable Matrix Functions, Math.
Research 37, Akademie, Berlin.

Moiseyev N.G.,  Popov, G.Ya. 1990. Exact solution of the problem of bending of a semi-infinite
plate completely bonded to an elastic half-space, Izv.AN SSSR, Solid Mechanics 25, 113-125.

Movchan, A. B., Willis, J. R. 1995. Dynamic weight functions for a moving crack. II. Shear loading. J. Mech. Phys. Solids 43, 1369-1383. 

Noble, B., 1988.  Methods Based on the Wiener-Hopf Technique,  Chelsea, New York.

Parton, V.Z.,  Perlin, P.I., 1982.  Integral Equations in Elasticity, Mir, Moscow.

Slepyan, L.I., 2002. Models and Phenomena in Fracture Mechanics, Springer, Berlin.

Vekua, N.P., 1967. Systems of Singular Integral Equations, Noordhoff, Groningen.

Willis, J.R., Movchan, A.B., 1995.   Dynamic weight functions for a moving crack. I. Mode I loading, J.
Mech. Phys. Solids 43,  319-341.

Zlatin, A.N.,  Khrapkov, A.A.   1986. A semi-infinite crack that is parallel to the boundary of an elastic half-plane,
 Dokl. Akad. Nauk SSSR 291, no.4, 810-813.

\end{document}